\documentclass{amsart}
\usepackage{amsmath,amsthm,amsfonts,amsbsy,verbatim}
\usepackage[all]{xy}

\theoremstyle{plain}
\newtheorem{lemma}{Lemma}[section]
\newtheorem{theorem}[lemma]{Theorem}

\numberwithin{equation}{section}

\newcommand{\qbin}[2]{\genfrac{[}{]}{0pt}{}{#1}{#2}}
\newcommand{\bin}[2]{\genfrac{(}{)}{0pt}{}{#1}{#2}}
\newcommand{\qbins}[2]{{\textstyle\genfrac{[}{]}{0pt}{}{#1}{#2}}}
\def\Z{\mathbb{Z}}
\def\Zp{\mathbb{Z}_+}
\def\bs{\boldsymbol}
\def\vm{\bs m}
\def\ve{{\bs e}}
\def\vn{\bs n}
\def\vet{\bs \eta}
\def\vmu{\bs \mu}
\def\I{\mathcal{I}}

\begin{document}

\title[Generalized $q$-Saalsch\"utz sum]{A generalization of the 
$q$-Saalsch\"utz sum and the Burge transform}

\author[A.~Schilling]{Anne Schilling}
\address{Instituut voor Theoretische Fysica, Universiteit van Amsterdam,
Valckenierstraat 65, 1018 XE Amsterdam, The Netherlands,
(Present address: Department of Mathematics,
Massachusetts Institute of Technology,
77 Massachusetts Ave, Cambridge MA 02139, USA}
\email{anne@math.mit.edu}

\author[S.~O.~Warnaar]{S.~Ole Warnaar}
\address{Instituut voor Theoretische Fysica, Universiteit van Amsterdam,
Valckenierstraat 65, 1018 XE Amsterdam, The Netherlands}
\email{warnaar@wins.uva.nl}
 


\begin{abstract}
A generalization of the $q$-(Pfaff)--Saalsch\"utz summation formula is proved.
This implies a generalization of the Burge transform, resulting in an 
additional dimension of the ``Burge tree''.
Limiting cases of our summation formula imply the (higher-level)
Bailey lemma, provide a new decomposition of the $q$-multinomial coefficients,
and can be used to prove the Lepowsky and Primc formula for the A$_1^{(1)}$
string functions.
\end{abstract}

\maketitle

\section{Introduction}
One of the most important summation formulas for basic hypergeometric
functions is Jackson's $q$-analogue of a $_3F_2$ summation formula
of Pfaff and Saalsch\"utz.
Employing standard notation (see e.g., Gaspar and Rahman~\cite{GR90})
this $q$-(Pfaff)--Saalsch\"utz sum is written as
\begin{equation}\label{qS1}
{_3\phi_2}\Bigl[\genfrac{}{}{0pt}{}{a,b,q^{-n}}{c,abq^{1-n}/c};q,q\Bigr]
:=\sum_{k=0}^n \frac{(a)_k(b)_k(q^{-n})_k\,q^k}{(q)_k(c)_k(abq^{1-n}/c)_k}
=\frac{(c/a)_n(c/b)_n}{(c)_n(c/ab)_n}
\end{equation}
for $n\in\Zp$.
Here $(a)_n$ is the $q$-shifted factorial, defined for all integers $n$ by
\begin{equation*}
(a;q)_{\infty}=(a)_\infty=\prod_{k=0}^{\infty}(1-aq^k) \qquad \text{and}
\qquad
(a;q)_n=(a)_n=\frac{(a)_{\infty}}{(aq^n)_{\infty}}.
\end{equation*}

Defining the $q$-binomial coefficient as
\begin{equation}\label{defqb}
\qbin{m+n}{m}=\begin{cases}\displaystyle \frac{(q)_{m+n}}{(q)_m(q)_n} &
\text{for $m,n\geq 0$}\\[2mm]
0 & \text{otherwise,}\end{cases}
\end{equation}
the $q$-Saalsch\"utz sum is often written as the following summation 
formula~\cite{Gould72,Carlitz74,Andrews76}
\begin{equation}\label{qS2}
\sum_{i=0}^M q^{i(i+\ell)}\qbin{L_1+L_2+M-i}{M-i}
\qbin{L_1}{i+\ell}\qbin{L_2}{i}=
\qbin{L_1+M}{M+\ell}\qbin{L_2+M+\ell}{M},
\end{equation}
valid for all $L_1,L_2,M,\ell\in\Z$ except when 
$-L_1\leq -\ell\leq L_2<0\leq M$ or 
$-L_2\leq\ell\leq L_1<0\leq M+\ell$. (In these cases the left-hand side
is zero whereas the right-hand side is not.)

In this paper we generalize the representation \eqref{qS2} of the 
$q$-Saalsch\"utz sum to a summation formula which transforms an $N$-fold
sum over a product of $N+2$ $q$-binomials to an $(N-1)$-fold sum over 
a product of $N+1$ $q$-binomials as stated in Theorem \ref{thm ex burge}
of the next section.
This generalized $q$-Saalsch\"utz sum contains many important special 
cases and can be applied in connection with the Burge transform,
the Bailey lemma, $q$-multinomial coefficients and level-$N$ A$_1^{(1)}$ 
string functions as summarized below.
\begin{enumerate}
\item
In Ref.~\cite{Burge93} Burge used equation~\eqref{qS2} to establish a
transformation on generating functions of (restricted) partition pairs. 
This ``Burge transform'', which generalizes a special case of the Bailey
lemma, can be used to derive a tree of identities for doubly bounded Virasoro 
characters~\cite{Burge93,FLW98}.
Our generalization of \eqref{qS2} adds a further dimension to the Burge 
tree as discussed in Section \ref{sec burge}.
\item
Letting $b$ tend to infinity in \eqref{qS1} yields the $q$-Chu--Vandermonde 
summation \cite[Eq. (II.7)]{GR90}. The $q$-binomial version of this is 
obtained by letting $M$ tend to infinity in \eqref{qS2}, resulting in
\begin{equation}\label{qCV}
\sum_{i=0}^{L_2}q^{i(i+\ell)}\qbin{L_1}{i+\ell}\qbin{L_2}{i}=
\qbin{L_1+L_2}{L_1-\ell}
\end{equation}
for $L_1,L_2,M,\ell\in\Z$ except when
$-L_1\leq -\ell\leq L_2<0$ or $-L_2\leq\ell\leq L_1<0$.
This identity can be viewed as a decomposition of the $q$-binomial
and is easily understood combinatorially using the notion of 
the Durfee rectangle of a partition.

The $q$-binomials have been generalized to $q$-trinomials in 
Ref.~\cite{AB87}, and more generally to $q$-multinomials in
Refs.~\cite{Andrews94,Butler94,Kirillov95,Schilling96,Warnaar97}.
Our generalized $q$-Saalsch\"utz sum implies a generalized 
$q$-Chu--Vandermonde sum which provides a new decomposition 
formula for $q$-multinomials in terms 
of $q$-binomials (see Section \ref{sec multinomials}).
\item
When $L_1$ and $L_2$ tend to infinity in \eqref{qS2} we are left with
\begin{equation}\label{cbpN1}
\sum_{i=0}^M \frac{q^{i(i+\ell)}}{(q)_{M-i}}
\frac{1}{(q)_i(q)_{i+\ell}}=\frac{1}{(q)_M(q)_{M+\ell}}.
\end{equation}
Let $\{\gamma\}_{L\ge 0}$ and $\{\delta\}_{L\ge 0}$ be sequences that satisfy
\begin{equation}\label{cBp}
\gamma_L=\sum_{r=L}^{\infty} \frac{\delta_r}{(q)_{r-L}(aq)_{r+L}}.
\end{equation}
Then the pair $(\gamma,\delta)$ is called a conjugate Bailey pair
relative to $a$~\cite{Bailey49,SW98}.
Replacing $M\to M-L$ and $\ell\to \ell+2L$ in equation \eqref{cbpN1} 
implies the conjugate Bailey pair
\begin{equation*}
\gamma_L=\frac{a^L q^{L^2}}{(q)_{M-L}(aq)_{M+L}}, \qquad
\delta_L=\frac{a^L q^{L^2}}{(q)_{M-L}},
\end{equation*}
with $a=q^{\ell}$.

A limit of our generalized $q$-Saalsch\"utz sum yields 
(a special case of) the higher-level generalization of this conjugate 
Bailey pair of Refs.~\cite{SW97,SW98}. For details see 
Section~\ref{sec bailey}.
This paper thus provides a new proof of the higher-level Bailey 
lemma of \cite{SW97,SW98} for a special choice of one of the parameters.

\item
Finally, letting $L_1, L_2$ and $M$ all tend to infinity in \eqref{qS2}
yields the well-known Durfee rectangle identity
\begin{equation*}
\sum_{i=0}^{\infty}\frac{q^{i(i+\ell)}}{(q)_i(q)_{i+\ell}}=
\frac{1}{(q)_{\infty}}.
\end{equation*}
This formula has many interpretations. Here we only mention that
the right-hand side can be identified with the level-1 A$_1^{(1)}$ 
string function.
Combined with the spinon formula of the string function
of Refs. \cite{ANOT96,BLS95,NY96a,NY96b,SW99},
the analogous limit of our generalized $q$-Saalsch\"utz sum yields the
fermionic expression for the string function due to Lepowsky and Primc
\cite{LP85} (see Section \ref{sec string}).
\end{enumerate}

\section{A generalized $q$-Saalsch\"utz identity}

The next theorem states the main result of this paper and provides
a generalization of the $q$-Saalsch\"utz summation formula \eqref{qS2}.
Let $C$ be the Cartan matrix of A$_{N-1}$ (i.e., $C_{ij}=2\delta_{i,j}-
\delta_{|i-j|,1}$ for $i,j=1,\dots,N-1$ where $\delta_{i,j}$ is the 
Kronecker delta symbol) and let $\I=2I-C$ be the 
corresponding incidence matrix where $I$ is the identity matrix.
Furthermore, let $\ve_i$, $i=1,\dots,N-1$ be the standard unit vectors in
$\Z^{N-1}$, $(\ve_i)_j=\delta_{i,j}$, and denote $\vn C^{-1} \vn
=\sum_{i,j=1}^{N-1}
n_i C^{-1}_{ij} n_j$ and $\ve_i C^{-1} \vn=(C^{-1} \vn)_i$
for $\vn\in\Z^{N-1}$.

\begin{theorem}\label{thm ex burge}
Let $\sigma=0,1$ and let $N,\ell,M,L_1+\frac{\ell+\sigma}{2},
L_2+\frac{\ell+\sigma}{2}$ be integers such that 
$\ell+\sigma N$ is even, $N\geq 1$ and $L_1,L_2\ge 0$. Then
\begin{multline}\label{sum}
\sum_{i=0}^M  q^{i(i+\ell)/N}  
\qbins{L_1+L_2+M-i}{M-i} 
\sum_{\substack{\vn\in \Z^{N-1} \\ 
\frac{2i+\ell+\sigma N}{2N}+(C^{-1} \vn)_1 \in \Z}}
q^{\vn C^{-1} \vn} \qbins{\vm+\vn}{\vn} 
 \qbins{L_1+\frac{1}{2}m_1}{i+\ell} 
 \qbins{L_2+\frac{1}{2}m_1}{i} \\[3mm]
=\sum_{\substack{\vet\in \Z^{N-1} \\ 
\frac{\ell+\sigma N}{2N}+(C^{-1} \vet)_1 \in \Z}}
q^{\vet C^{-1} \vet} \qbins{\vmu+\vet}{\vet}
 \qbins{L_1+\frac{1}{2}(M+\mu_1)}{M+\ell}
 \qbins{L_2+\frac{1}{2}(M+\ell+\mu_{N-1})}{M}
\end{multline}
with 
\begin{equation}\label{mn}
\vm+\vn=\frac{1}{2}(\I \vm+(2i+\ell)\ve_1)
\end{equation}
and
\begin{equation}\label{etamu}
\vmu+\vet=\frac{1}{2}(\I \vmu+(M+\ell)\ve_1+M\ve_{N-1}).
\end{equation}
\end{theorem}

The vector $\vm\in\Z^{N-1}$ on the left-hand side 
is determined by the (summation) variable $\vn$ through 
the $(\vm,\vn)$-system \eqref{mn}. 
Similarly $\vmu\in\Z^{N-1}$ is determined by \eqref{etamu}.
Also, $\qbins{\vm+\vn}{\vn}=\prod_{j=1}^{N-1}
\qbins{m_j+n_j}{n_j}$ and similarly for $\qbins{\vmu+\vet}{\vet}$.
We further note that the nature of the solutions of \eqref{mn}
depends on the parity of $N$. When $N$ is odd one must have
$$ m_1\equiv m_3\equiv \dots\equiv m_{N-2}\equiv 0 \pmod{2},\quad
m_2\equiv m_4\equiv \dots\equiv m_{N-1}\equiv \ell \pmod{2}$$
whereas for $N$ even one finds
\begin{equation}\label{parities}
m_1\equiv m_3\equiv \dots\equiv m_{N-1} \pmod{2}, \quad 
m_2\equiv m_4\equiv \dots\equiv m_{N-2}\equiv \ell\equiv 0 \pmod{2}.
\end{equation}
This implies that $m_1$ is even for $N$ odd so that $L_1,L_2$ must
be integers. This indeed follows from (since $N$ is odd) $0\equiv\ell+\sigma N
\equiv \ell+\sigma\pmod{2}$. 
When $N$ is even the partity of $m_1$ is not fixed and
there is the freedom to choose $m_1$ even corresponding to $\sigma=0$
or $m_1$ odd corresponding to $\sigma=1$.
(Since for $N$ even $0\equiv \ell+\sigma N\equiv \ell\pmod{2}$, $\ell$ 
is even in accordance with \eqref{parities} and hence, since $L_1,L_2$ 
must be integers when $m_1$ even and half an odd integer when $m_1$ odd,
it thus follows from $L_i+(\ell+\sigma)/2\in \Z$
that $\sigma$ has the same parity as $m_1$.)
A similar analysis of the solutions of the $(\vmu,\vet)$-system \eqref{etamu}
can be carried out.
The restrictions on the sums over $\vn$ and $\vet$ ensure that
the components of $\vm$ and $\vmu$ are integer and have the parity as 
discussed above. 

Equation \eqref{sum} yields a summation formula for every $N\ge 1$.
When $N=1$ the sums over $\vn$ and $\vet$ drop out; on the left-hand
side $m_1=0$ and on the right-hand side one needs to interpret $\mu_1=M$
and $\mu_0=M+\ell$. Then \eqref{sum} indeed reduces to \eqref{qS2} 
for $N=1$.

\begin{proof}[Proof of Theorem \ref{thm ex burge}]
Note that both sides of \eqref{sum} are zero
unless $M+\ell\ge 0$ and $M\ge 0$. Furthermore, denoting the identity 
\eqref{sum} by $I(L_1,L_2,M,\ell)$, it enjoys the symmetry
$I(L_1,L_2,M,\ell)=I(L_2,L_1,M+\ell,-\ell)$. Hence we may assume 
$\ell\ge 0$ and $M\ge 0$ in the proof below.

Throughout the proof we use modified $q$-binomials defined as
\begin{equation}\label{mod qbin}
\qbin{m+n}{m}=\frac{(q^{n+1})_m}{(q)_m} \qquad 
\text{for $m\in\Zp$, $n\in\Z$,}
\end{equation}
and zero otherwise. Note that $\qbins{m+n}{m}$ is zero if $n<0$ unless 
$m+n<0$.
Let us now show that on both sides of \eqref{sum} 
the $q$-binomials \eqref{defqb}
can be replaced by the modified $q$-binomials. Since $M,\ell,L_1,L_2\ge 0$
we find from \eqref{mn} and \eqref{etamu} that $m_i+n_i\ge 0$
and $\mu_i+\eta_i\ge 0$ if $m_i,\mu_i\ge 0$ so that $\qbins{\vm+\vn}{\vn}$ 
and $\qbins{\vmu+\vet}{\vet}$ in \eqref{sum} can be replaced by
the modified $q$-binomials $\qbins{\vm+\vn}{\vm}$ and
$\qbins{\vmu+\vet}{\vmu}$, respectively. The other $q$-binomials can
be turned into modified $q$-binomials since the top entries are nonnegative 
by the conditions on the parameters.

The proof of \eqref{sum} makes frequent use of the following identity which
is a corollary of Sears' transformation formula for a balanced $_4\phi_3$
series \cite[Eq. (III.15)]{GR90}
\begin{multline}\label{sears}
\sum_{i\in\Z}q^{i(i-a+e+g)}\qbins{i+a}{a}\qbins{b-i}{c-i}\qbins{d}{i+e}
\qbins{f}{i+g}\\
=\sum_{i\in\Z}q^{i(i-a+e+g)}\qbins{a-g}{a-g-i}\qbins{b-d+e}{c-i}
\qbins{c+d-i}{c+e}\qbins{i+f}{i+g},
\end{multline}
where $a,b,c,d,e,f,g\in\Z$ and the condition $a+b=c+d+f$ applies. 
Since we need the Sears transform \eqref{sears} with negative entries 
in the $q$-binomials it is
essential that definition \eqref{mod qbin} is used here.
(The above formula is not correct for all $a,\dots,g\in\Z$ with the
use of \eqref{defqb}).

We start by shifting $\vn\to \vn+i \ve_1$, followed by $i\to i-n_1$.
This transforms the left-hand side of \eqref{sum} into
\begin{multline*}
\sum_{i,\vn} q^{(i-n_1)(i-n_1-m_1+\ell)+\vn C^{-1} \vn}\\ \times 
\qbins{L_1+L_2+M+n_1-i}{M+n_1-i} 
\qbins{L_1+\frac{1}{2}m_1}{i+\ell-n_1}
\qbins{L_2+\frac{1}{2}m_1}{i-n_1}
\qbins{m_1+i}{m_1} 
\prod_{\alpha=2}^{N-1}\qbins{m_{\alpha}+n_{\alpha}}{m_{\alpha}} 
\end{multline*}
where the sum over $\vn$ is restricted by
\begin{equation}\label{res} 
\frac{\ell+\sigma N}{2N}+(C^{-1}\vn)_1\in\Z
\end{equation}
and the $(\vm,\vn)$-system is given by
\begin{equation}\label{mnnew}
\vm+\vn=\frac{1}{2}(\I \vm+\ell \ve_1).
\end{equation}
Since the $(\vm,\vn)$-system has become $i$-independent, only the
first four $q$-binomials depend on the summation variable $i$.
Hence we may apply \eqref{sears} with
$a=m_1$, $b=L_1+L_2+M+n_1$, $c=M+n_1$, $d=L_1+\frac{1}{2}m_1$,
$e=\ell-n_1$, $f=L_2+\frac{1}{2}m_1$ and $g=-n_1$ to obtain
\begin{multline*}
\sum_{i,\vn} q^{(i-n_1)(i-n_1-m_1+\ell)+\vn C^{-1} \vn}\\ \times
\qbins{L_1+M+n_1+\frac{1}{2}m_1-i}{M+\ell}
\qbins{L_2+M+\ell-\frac{1}{2}m_1}{M+n_1-i}
\qbins{L_2+\frac{1}{2}m_1+i}{i-n_1}
\qbins{m_1+n_1}{m_1+n_1-i}
\prod_{\alpha=2}^{N-1}\qbins{m_{\alpha}+n_{\alpha}}{m_{\alpha}} .
\end{multline*}
Shifting $\vn\to \vn+i(2\ve_1-\ve_2)$ and $\vm\to \vm-2i \ve_1$, 
which leaves the $(\vm,\vn)$-system \eqref{mnnew} and the restriction 
\eqref{res} on the summation over $\vn$ invariant, yields
\begin{multline*}
\sum_{i,\vn} q^{(i+\frac{m_2-m_1}{2})^2-(\frac{m_1-\ell}{2})^2+\vn C^{-1} \vn}
\qbins{L_1+M+\frac{1}{2}m_1+n_1}{M+\ell}
\qbins{L_2+M+\ell-\frac{1}{2}m_1+i}{M+n_1+i} \\
\times \qbins{L_2+\frac{1}{2}m_1}{-i-n_1}
\qbins{m_1+n_1}{m_1+n_1-i}\qbins{m_2+n_2-i}{m_2} 
\prod_{\alpha=3}^{N-1}\qbins{m_{\alpha}+n_{\alpha}}{m_{\alpha}},
\end{multline*}
where we have used the $(\vm,\vn)$-system to simplify the exponent of $q$.
Shifting $i\to n_2-i$ one can apply \eqref{sears}
with $a=m_2$, $b=L_2+M+\ell-\frac{1}{2}m_1+n_2$,
$c=M+n_1+n_2$, $d=L_2+\frac{1}{2}m_1$, $e=-n_1-n_2$,
$f=m_1+n_1$ and $g=m_1+n_1-n_2$,
observing that $$c+d+f-a-b=2m_1+2n_1-m_2-\ell=0$$ thanks to \eqref{mnnew}.
This yields
\begin{multline*}
\sum_{i,\vn} q^{(i-\frac{m_3-m_2}{2})^2-(\frac{m_1-\ell}{2})^2+\vn C^{-1} \vn}
\qbins{L_1+M+\frac{1}{2}m_1+n_1}{M+\ell}
\qbins{L_2+M+\frac{1}{2}m_1+n_1+n_2-i}{M}\\
\times\qbins{m_2+n_2-m_1-n_1}{m_2+n_2-m_1-n_1-i}
\qbins{M+\ell-m_1-n_1}{M+n_1+n_2-i}
\qbins{m_1+n_1+i}{m_1+n_1-n_2+i}
\prod_{\alpha=3}^{N-1}\qbins{m_{\alpha}+n_{\alpha}}{m_{\alpha}}.
\end{multline*}
Shifting $\vn\to \vn+i(\ve_1+\ve_2-\ve_3)$ and 
$\vm\to \vm-2i(\ve_1+\ve_2)$, which again leaves the $(\vm,\vn)$-system 
\eqref{mnnew} and the restriction \eqref{res}
on the sum over $\vn$ unchanged, leads to
\begin{multline}\label{f3}
\sum_{i,\vn} q^{(i+\frac{m_3-m_2}{2})^2-(\frac{m_1-\ell}{2})^2+\vn C^{-1} \vn}
\qbins{L_1+M+n_1+\frac{1}{2}m_1}{M+\ell}
\qbins{L_2+M+\frac{1}{2}m_1+n_1+n_2}{M} \\
\times
\qbins{m_2+n_2-m_1-n_1}{m_2+n_2-m_1-n_1-i}
\qbins{M+\ell-m_1-n_1+i}{M+n_1+n_2+i}
\qbins{m_1+n_1}{m_1+n_1-n_2-i}
\qbins{m_3+n_3-i}{m_3}
\prod_{\alpha=4}^{N-1}\qbins{m_{\alpha}+n_{\alpha}}{m_{\alpha}}.
\end{multline}
We now need the following lemma.
\begin{lemma}\label{lem help}
For $p=3,\dots,N$, let
\begin{equation*}
\begin{split}
f_p=&\sum_{i,\vn}  q^{(i+\frac{m_p-m_{p-1}}{2})^2-
(\frac{m_1-\ell}{2})^2+\vn C^{-1} \vn} 
\qbins{L_1+M+\frac{1}{2}m_1+n_1}{M+\ell}
\qbins{L_2+M+\frac{1}{2}m_1+n_1+n_2}{M} \\
&\times \Bigl(\prod_{\alpha=1}^{p-3} 
\qbins{M+\sum_{\beta=1}^{\alpha+2}n_{\beta}
+\sum_{\beta=1}^{\alpha}(-1)^{\alpha-\beta}
(m_\beta+n_\beta)}{M+\sum_{\beta=1}^{\alpha}n_{\beta}
+\sum_{\beta=1}^{\alpha}(-1)^{\alpha-\beta}
(m_\beta+n_\beta)}\Bigr) \Bigl(
\prod_{\alpha=p+1}^{N-1}\qbins{m_{\alpha}+n_{\alpha}}{m_{\alpha}} \Bigr) \\
& \times 
\qbins{\sum_{\alpha=1}^{p-1}(-1)^{p-\alpha-1}(m_\alpha+n_\alpha)}
{\sum_{\alpha=1}^{p-1}(-1)^{p-\alpha-1}(m_\alpha+n_\alpha)-i}
\qbins{M+\ell-m_1-\sum_{\alpha=1}^{p-2}n_{\alpha}+i}
{M+\sum_{\alpha=1}^{p-1} n_{\alpha}+i}\\[2mm]
&\times
\qbins{\sum_{\alpha=1}^{p-2}(-1)^{p-\alpha}(m_\alpha+n_\alpha)}
{\sum_{\alpha=1}^{p-2}(-1)^{p-\alpha}(m_\alpha+n_\alpha)-n_{p-1}-i}
\qbins{m_p+n_p-i}{m_p},
\end{split}
\end{equation*}
with $(\vm,\vn)$-system \eqref{mnnew} and $m_N=n_N=0$.
Then $f_p=f_{p+1}$ for $3\le p<N$.
\end{lemma}
\begin{proof}
Change $i\to n_p-i$ and apply \eqref{sears} with
$a=m_p$, $b=M+\ell-m_1-\sum_{\alpha=1}^{p-2}n_\alpha+n_p$,
$c=M+\sum_{\alpha=1}^p n_\alpha$, $d=\sum_{\alpha=1}^{p-2}(-1)^{p-\alpha}
(m_\alpha+n_\alpha)$, $e=d-n_{p-1}-n_p$,
$f=\sum_{\alpha=1}^{p-1}(-1)^{p-\alpha-1}(m_\alpha+n_\alpha)$
and $g=f-n_p$,
observing that
$$c+d+f-a-b=2\sum_{\alpha=1}^{p-1}n_{\alpha}+m_1+m_{p-1}-m_p-\ell=0$$
by summing up the first $p-1$ components of the $(\vm,\vn)$-system
\eqref{mnnew}.
This leads to
\begin{equation}\label{fp1}
\begin{split}
f_p=&\sum_{i,\vn}  q^{(n_p-i+\frac{m_p-m_{p-1}}{2})^2-
(\frac{m_1-\ell}{2})^2+\vn C^{-1} \vn} 
\qbins{L_1+M+\frac{1}{2}m_1+n_1}{M+\ell}
\qbins{L_2+M+\frac{1}{2}m_1+n_1+n_2}{M} \\
&\times
\Bigl(\prod_{\alpha=1}^{p-3} 
\qbins{M+\sum_{\beta=1}^{\alpha+2}n_{\beta}
+\sum_{\beta=1}^{\alpha}(-1)^{\alpha-\beta}(m_\beta+n_\beta)}
{M+\sum_{\beta=1}^{\alpha}n_{\beta}
+\sum_{\beta=1}^{\alpha}(-1)^{\alpha-\beta}(m_\beta+n_\beta)}\Bigr) 
\Bigl(\prod_{\alpha=p+1}^{N-1}
\qbins{m_{\alpha}+n_{\alpha}}{m_{\alpha}} \Bigr)\\
&\times
\qbins{\sum_{\alpha=1}^p(-1)^{p-\alpha}(m_\alpha+n_\alpha)}
{\sum_{\alpha=1}^p(-1)^{p-\alpha}(m_\alpha+n_\alpha)-i}
\qbins{M+\ell-m_1-\sum_{\alpha=1}^{p-1}n_{\alpha}}
{M+\sum_{\alpha=1}^p n_{\alpha}-i}\\[2mm]
&\times
\qbins{M+\sum_{\alpha=1}^p n_{\alpha}+\sum_{\alpha=1}^{p-2}(-1)^{p-\alpha}
(m_\alpha+n_\alpha)-i}
{M+\sum_{\alpha=1}^{p-2} n_{\alpha}+\sum_{\alpha=1}^{p-2}(-1)^{p-\alpha}
(m_\alpha+n_\alpha)}
\qbins{\sum_{\alpha=1}^{p-1}(-1)^{p-\alpha-1}(m_\alpha+n_\alpha)+i}
{\sum_{\alpha=1}^{p-1}(-1)^{p-\alpha-1}(m_\alpha+n_\alpha)-n_p+i}.
\end{split}
\end{equation}
We now carry out the transformations
$\vn\to \vn+i(\ve_1+\ve_p-\ve_{p+1})$ and 
$\vm\to \vm-2i(\ve_1+\ve_2+\cdots+\ve_p)$, which leave the 
$(\vm,\vn)$-system unchanged. (Here $\ve_N:=0$.)
Using $\vn C^{-1}(\ve_1+\ve_p-\ve_{p+1})=\sum_{\alpha=1}^p n_{\alpha}$
and $(\ve_1+\ve_p-\ve_{p+1})C^{-1}(\ve_1+\ve_p-\ve_{p+1})=2$, as well as the
$(\vm,\vn)$-system, yields
\begin{multline*}
(n_p-i+\frac{m_p-m_{p-1}}{2})^2-
(\frac{m_1-\ell}{2})^2+\vn C^{-1} \vn \to \\
(n_p+\frac{m_p-m_{p-1}}{2})^2-
(i-\frac{m_1-\ell}{2})^2+2i(i+\sum_{\alpha=1}^p n_{\alpha})+\vn C^{-1} \vn\\
=(i+\frac{m_{p+1}-m_p}{2})^2-(\frac{m_1-\ell}{2})^2+\vn C^{-1} \vn
\end{multline*}
transforming \eqref{fp1} into $f_{p+1}$ as desired.
\end{proof}
Equation \eqref{f3} corresponds to $f_3$ and we can thus use the above lemma
to replace it with $f_N$. Since $m_N=0$, the last $q$-binomial in $f_N$ is 1
and we can perform the sum over $i$ using the $q$-Saalsch\"utz sum, 
which is the special case $a=0$ of the Sears transformation \eqref{sears}.
(When $a=0$, the only nonvanishing term on the right-hand side of
\eqref{sears} corresponds to $i=-g$.)
Specifically, we take $f_N$, replace $i$ by $-i$ and apply \eqref{sears}
with the same choice of parameters as in the proof of Lemma~\ref{lem help}
but with $p=N$, $n_N=0$ and $a=m_N=0$. Then we get
\begin{equation}\label{mess}
\begin{split}
&\sum_{\vn}  q^{(\frac{m_{N-1}}{2})^2-
(\frac{m_1-\ell}{2})^2
+(n_{N-1}-\sum_{\alpha=1}^{N-2}(-1)^{N-\alpha}(m_\alpha+n_\alpha))
(\sum_{\alpha=1}^{N-1}(-1)^{N-\alpha-1}(m_\alpha+n_\alpha))}\\
&\qquad\times q^{\vn C^{-1} \vn}
\Bigl(\prod_{\alpha=1}^{N-3} 
 \qbins{M+\sum_{\beta=1}^{\alpha+2}n_{\beta}
+\sum_{\beta=1}^{\alpha}(-1)^{\alpha-\beta}(m_\beta+n_\beta)}
{M+\sum_{\beta=1}^{\alpha}n_{\beta}
+\sum_{\beta=1}^{\alpha}(-1)^{\alpha-\beta}(m_\beta+n_\beta)}\Bigr)\\
&\qquad\times 
 \qbins{L_1+M+\frac{1}{2}m_1+n_1}{M+\ell}
 \qbins{M+\ell-m_1-\sum_{\alpha=1}^{N-2}n_{\alpha}}
  {M+\ell-m_1-\sum_{\alpha=1}^{N-1}n_{\alpha}
  -\sum_{\alpha=1}^{N-1}(-1)^{N-\alpha-1}(m_\alpha+n_\alpha)}\\[2mm]
&\qquad\times
 \qbins{L_2+M+\frac{1}{2}m_1+n_1+n_2}{M}
 \qbins{M+\ell-m_1-\sum_{\alpha=1}^{N-1}n_{\alpha}}
  {M+\ell-m_1-\sum_{\alpha=1}^{N-2}n_{\alpha}
  -\sum_{\alpha=1}^{N-2}(-1)^{N-\alpha}(m_\alpha+n_\alpha)}.
\end{split}
\end{equation}
All that remains to be done is to clean up the above expression.
Introduce a new variable $\vet\in\Z^{N-1}$ through its components
as follows
\begin{align*}
\eta_i&=n_{2i}+n_{2i+1}  &&\text{for $i=1,\dots,\lfloor N/2 \rfloor -1$}\\
\eta_{N-i}&=n_{2i+1}+n_{2i+2} &&
 \text{for $i=1,\dots,\lfloor (N-1)/2 \rfloor -1$} \\
\eta_{\lfloor (N+1)/2 \rfloor}&=
\sum_{\alpha=1}^{N-2}(-1)^{N-\alpha}(m_\alpha+n_\alpha)-n_{N-1} &\\
\eta_{\lfloor (N+1)/2 \rfloor \pm 1}&=
\sum_{\alpha=1}^{N-1}(-1)^{N-\alpha-1}(m_\alpha+n_\alpha)+n_{N-1}&
\end{align*}
for $N$ even/odd.
Also define $\vmu$ through the $(\vmu,\vet)$-system \eqref{etamu}
Eliminating $\vm$ and $\vn$ from \eqref{mess} in favour of $\vmu$ and $\vet$,
we finally get the right-hand side of \eqref{sum}.
We also note that $(C^{-1}\vet)_1=\sum_{i=1}^{N-1} (N-i)\eta_i/N$ yields
$(-n_1+\sum_{i=2}^{N-1} (N-i)n_i)/N=(C^{-1}\vn)_1-n_1$ so that
the restriction \eqref{res} on the sum over $\vn$ translates
into the restriction
\begin{equation*}
\frac{\ell+\sigma N}{2N}+(C^{-1}\vet)_1\in\Z
\end{equation*}
for the sum over $\vet$ as it should.

\end{proof}

\section{The Burge transform}\label{sec burge}

Perhaps the most interesting application of our generalized
$q$-Saalsch\"utz sum \eqref{sum} arises when it is combined with
the Burge transform~\cite{Burge93,FLW98}.
The Burge transform is a generalization of (a special case)
of the Bailey lemma and can be utilized to derive an
infinite tree (a Burge tree) of polynomial identities from
a single initial identity. In this section we show that
each element of a Burge tree can be transformed using \eqref{sum}
to yield an additional infinite series of polynomial identities.

In his study of restricted partition pairs Burge considered the polynomial
\begin{multline}\label{X}
X_{r,s}^{(p,p')}(M_1,L_1,M_2,L_2) \\
\begin{split} =\sum_{j=-\infty}^{\infty}
\Bigl\{& q^{j(pp'j+p'(M_{12}+r)-ps)}
\qbins{M_1+L_1-(p'-p)j}{M_1+pj}\qbins{M_2+L_2+(p'-p)j}{M_2-pj} \\
-&q^{(pj+M_{12}+r)(p'j+s)}
\qbins{M_1+L_1-(p'-p)j+r-s}{M_1+pj+r}
\qbins{M_2+L_2+(p'-p)j-r+s}{M_2-pj-r}\Bigr\},
\end{split}
\end{multline}
with $M_{12}=M_1-M_2$,
and proved that it is the generating function of pairs of partitions
$(\lambda,\mu)$ such that
\begin{equation*}
0\leq \lambda_1\leq \dots\leq \lambda_{M_1}\leq L_1, \qquad
0\leq \mu_1\leq \dots\leq \mu_{M_2}\leq L_2, \qquad
\end{equation*}
and
\begin{equation*}
\lambda_i-\mu_{i-r+1}\geq 1-s,\qquad \mu_i-\lambda_{i-p+r+1}\geq 1-p'+s.
\end{equation*}
Here the integers $p,p',r,s$ are restricted to $p,p'\geq 1$, 
$0\leq r+M_{12}\leq p$ and $0\leq s-L_{12}\leq p'$, with $L_{12}=L_1-L_2$.
There are four exceptional cases, $r=0$, $r=p$, $r=-M_{12}$ and $r=p-M_{12}$
that demand the additional conditions $\mu_1\leq s-1$, $\lambda_1\leq p'-s-1$,
$\lambda_{M_2}\geq L_1-s+1$ and $\mu_{M_1}\geq L_2-p'+s+1$, 
respectively~\cite{GK97,FLW98}.

The important observation made in \cite{Burge93} is that
\begin{multline}\label{Bt}
X_{r,r+s}^{(p,p+p')}(M_1,L_1,M_2,L_2) \\
=\sum_{i\in\Z} q^{i(i+M_{12})}
\qbin{L_1+L_2+M_2-i}{M_2-i} X_{r,s}^{(p,p')}(i+M_{12},L_1-i,i,L_2-M_{12}-i)
\end{multline}
and
\begin{multline}\label{Bt2}
X_{s-M_{12},r+s+L_{12}}^{(p',p+p')}(M_1,L_1,M_2,L_2) \\
=\sum_{i\in\Z} q^{i(i+M_{12})}
\qbin{L_1+L_2+M_2-i}{M_2-i} 
X_{r,s}^{(p,p')}(L_1-i,i+M_{12},L_2-M_{12}-i,i)
\end{multline}
where the second equation follows from the first by exploiting
the symmetry
\begin{equation}\label{symmetry}
X_{r,s}^{(p,p')}(M_1,L_1,M_2,L_2)=
X_{s-L_{12},r+M_{12}}^{(p',p)}(L_1,M_1,L_2,M_2).
\end{equation}
The proof of the 
Burge transform follows from the $q$-Saalsch\"utz formula \eqref{qS2}.
In \cite{Burge93,FLW98} the defining equation \eqref{X} is
substituted into \eqref{Bt}, then the sums over $i$ and $j$ are interchanged,
followed by the variable change $i\to i+pj$ and $i\to i+pj+r$ in the
terms corresponding to the second and third line of \eqref{X}, respectively
(referred to as the positive and negative terms below).
Then the $q$-Saalsch\"utz sum is used with $L_1\to L_1+M_{12}-(p'-p)j$,
$L_2\to L_2-M_{12}+(p'-p)j$, $M\to M_2-pj$ and $\ell\to M_{12}+2pj$
for the positive terms and
$L_1\to L_1+M_{12}-(p'-p)j+r-s$,
$L_2\to L_2-M_{12}+(p'-p)j-r+s$, $M\to M_2-pj-r$ and $\ell\to M_{12}+2pj+2r$
for the negative terms.
This gives the left-hand side of \eqref{Bt}.
However, we note that it needs to be verified that the summation
\eqref{qS2} has not been employed when the
variables therein lie in the ranges given just below \eqref{qS2}.
This means that
\begin{multline}\label{c1}
-L_1-M_{12}+(p'-p)j-r+s\leq -M_{12}-2pj-2r  \\
\leq L_2-M_{12}+(p'-p)j-r+s<0\leq M_2-pj-r
\end{multline}
and
\begin{multline}\label{c2}
-L_2+M_{12}-(p'-p)j+r-s\leq M_{12}+2pj+2r \\
\leq L_1+M_{12}-(p'-p)j+r-s<0\leq M_1+pj+r
\end{multline}
and the corresponding inequalities obtained by setting $r=s=0$ should
not hold for any $j\in\Z$.
Eliminating $j$ gives several conditions on the parameters in \eqref{Bt}.
In particular \eqref{c1} can only hold if 
$$ 2pj>-M_{12}-2r \quad \text{and} \quad 2p'j<M_{12}+L_{12}-2s. $$
Similarly, \eqref{c2} can only hold if
$$ 2pj<-M_{12}-2r \quad \text{and} \quad 2p'j>M_{12}+L_{12}-2s. $$
If, for example, $M_{12}=L_{12}=0$ these conditions cannot be satisfied for 
any $j$ recalling that $0\leq r\leq p$ and $0\leq s\leq p'$. 
Hence, setting 
$$X_{r,s}^{(p,p')}(M,L,M,L)=X_{r,s}^{(p,p')}(M,L),$$
the symmetric version of the Burge transform \eqref{Bt}
\begin{equation}
X_{r,r+s}^{(p,p+p')}(M,L)=\sum_{i=0}^M q^{i^2}
\qbin{2L+M-i}{2L} X_{r,s}^{(p,p')}(i,L-i)
\end{equation}
always holds.
By the same arguments one can show that the symmetric form of \eqref{Bt2}
\begin{equation*}
X_{s,r+s}^{(p',p+p')}(M,L)=\sum_{i=0}^M q^{i^2}
\qbin{2L+M-i}{2L} X_{r,s}^{(p,p')}(L-i,i)
\end{equation*}
is true for arbitrary $M$ and $L$.

By iterating the two Burge transformations, starting with an appropriate
initial identity for $X_{r,s}^{(p,p')}$, one can derive an infinite
tree of polynomial identities. This was mentioned in \cite{Burge93}
and explicitly carried out in \cite{FLW98}.
To illustrate this we follow \cite{FLW98} and use the trivial result
\begin{equation}\label{initial}
X_{0,1}^{(1,2)}(M,L)=\delta_{L,0}
\end{equation}
to derive the Burge tree
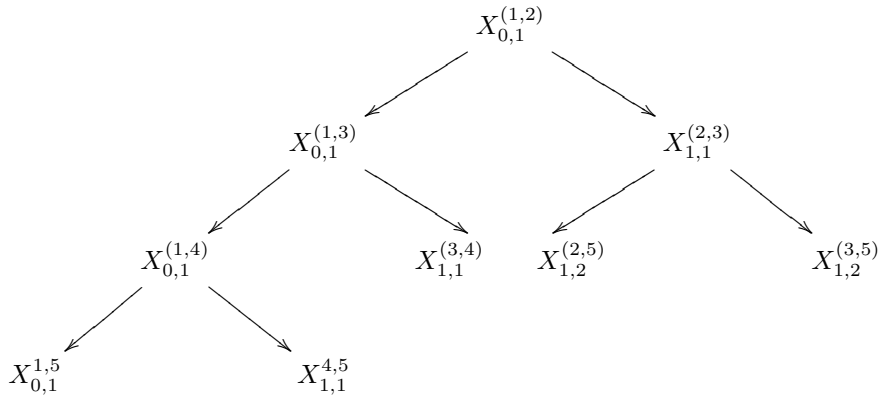
\begin{equation*}
\quad\xymatrix{
&&& X_{0,1}^{(1,2)} \ar[dl] \ar[dr]&&&\\
&& X_{0,1}^{(1,3)} \ar[dl] \ar[dr] &&
X_{1,1}^{(2,3)} \ar[dl] \ar[dr] &&\\
&X_{0,1}^{(1,4)} \ar[dr]\ar[dl] && \hspace{-3mm}X_{1,1}^{(3,4)} \qquad
X_{1,2}^{(2,5)}\hspace{-3mm} && X_{1,2}^{(3,5)} & \\
X_{0,1}^{1,5}&&X_{1,1}^{4,5}&&&&}
\end{equation*}
where the node labeled $X_{r,s}^{(p,p')}$ denotes a polynomial identity for
$X_{r,s}^{(p,p')}(M,L)$.
(Actually, in Ref.~\cite{FLW98} an extension of the Burge tree was
constructed by exploiting various symmetries of $X_{r,s}^{(p,p')}$.)
Explicitly some of the identities in the above tree are~\cite{Burge93,FLW98},
\begin{align}
\label{nn}
X_{0,1}^{(1,3)}(M,L)&=q^{L^2}\qbin{L+M}{2L} \\[1mm]
\label{Euler}
X_{1,1}^{(2,3)}(M,L)&=\qbin{2L+M}{2L} \\
\label{Ising}
X_{1,1}^{(3,4)}(M,L)&=\sum_{\substack{m=0\\ m~\text{even}}}^L
q^{\frac{1}{2}m^2}\qbin{2L+M-\frac{1}{2}m}{2L}\qbin{L}{m} \\
\label{RR}
X_{1,2}^{(2,5)}(M,L)&=\sum_{n=0}^Lq^{n^2}\qbin{2L+M-n}{2L}\qbin{2L-n}{n}.
\end{align}
Equation \eqref{Euler} is a doubly bounded version of the Euler identity,
equation \eqref{Ising} is a doubly bounded analogue of
the vacuum-character identity of the Ising model
$$
\sum_{\substack{m=0\\ m~\text{even}}}^{\infty}
\frac{q^{\frac{1}{2}m^2}}{(q)_m}=\frac{1}{2}
\Bigl\{(-q^{1/2})_{\infty}+(q^{1/2})_{\infty}\Bigr\}=
\prod_{j=1}^{\infty}\frac{(1+q^{8j-3})(1+q^{8j-5})(1-q^{8j})}
{1-q^{2j}}
$$
and \eqref{RR} is a doubly bounded version of the 
(first) Rogers--Ramanujan identity
$$
\sum_{n=0}^{\infty}\frac{q^{n^2}}{(q)_n}=
\prod_{j=1}^{\infty}\frac{1}{(1-q^{5j-1})(1-q^{5j-4})}.
$$

To see how \eqref{sum} transforms an identity in the Burge tree,
let us first introduce a generalization of the polynomial
$X_{r,s}^{(p,p')}(M_1,L_1,M_2,L_2)$ as follows. Let $N$ be a positive integer,
$\sigma=0,1$ and let $M_1,M_2,L_1+\frac{M_{12}+\sigma}{2},L_2
+\frac{M_{12}+\sigma}{2}$ be integers such that
$M_{12}+\sigma N$ is even. 
Also assume that $(p'-p)/N\in\Zp$ and $(r-s)/N\in\Z$, for $r,s$ integers.
Then
\begin{comment}
\begin{eqnarray*}
\lefteqn{
X_{r,s,\sigma}^{(p,p'),N}(M_1,L_1,M_2,L_2)} \\
& &=\sum_{j=-\infty}^{\infty}
q^{\frac{j}{N}(pp'j+p'(M_{12}+r)-ps)} 
\sum_{\substack{\vet\in \Z^{N-1} \\
\frac{M_{12}+2pj+\sigma N}{2N}+(C^{-1} \vet)_1 \in \Z}}
q^{\vet C^{-1} \vet} \qbins{\vet+\vmu}{\vet} \\
& & \times
\qbins{M_1+L_1+\frac{N-1}{2N}M_{12}-(p'/N-p)j-(C^{-1}\vet)_1}{M_1+pj}
\qbins{M_2+L_2-\frac{N-1}{2N}M_{12}+(p'/N-p)j-(C^{-1}\vet)_{N-1}}{M_2-pj} \\
& &-\sum_{j=-\infty}^{\infty}
q^{\frac{1}{N}(pj+M_{12}+r)(p'j+s)} 
\sum_{\substack{\vet\in \Z^{N-1} \\
\frac{M_{12}+2pj+2r+\sigma N}{2N}+(C^{-1} \vet)_1 \in \Z}}
q^{\vet C^{-1} \vet} \qbins{\vet+\vmu}{\vet} \\
& &\times 
\qbins{M_1+L_1+\frac{N-1}{2N}M_{12}-(p'/N-p)j+r-s/N-(C^{-1}\vet)_1}{M_1+pj+r}
\qbins{M_2+L_2-\frac{N-1}{2N}M_{12}+(p'/N-p)j-r+s/N-(C^{-1}\vet)_{N-1}}
{M_2-pj-r},
\end{eqnarray*}
\end{comment}
\begin{align}\label{XN}
X_{r,s,\sigma}^{(p,p'),N}&(M_1,L_1,M_2,L_2) \\
&=\sum_{j=-\infty}^{\infty}
q^{\frac{j}{N}(pp'j+p'(M_{12}+r)-ps)}
\sum_{\substack{\vet\in \Z^{N-1} \notag \\
\frac{M_{12}+2pj+\sigma N}{2N}+(C^{-1} \vet)_1 \in \Z}}
q^{\vet C^{-1} \vet} \qbins{\vet+\vmu}{\vet}\notag \\
&  \quad \times
\qbins{M_1+L_1-(p'-p)j/N-\frac{1}{2}(M_2-pj-\mu_1)}{M_1+pj}
\qbins{M_2+L_2+(p'-p)j/N-\frac{1}{2}(M_1+pj-\mu_{N-1})}{M_2-pj} \notag \\
& -\sum_{j=-\infty}^{\infty}
q^{\frac{1}{N}(pj+M_{12}+r)(p'j+s)}
\sum_{\substack{\vet\in \Z^{N-1} \\
\frac{M_{12}+2pj+2r+\sigma N}{2N}+(C^{-1} \vet)_1 \in \Z}}
q^{\vet C^{-1} \vet} \qbins{\vet+\vmu}{\vet}  \notag \\
& \quad \times
\qbins{M_1+L_1-((p'-p)j-r+s)/N-\frac{1}{2}(M_2-pj-r-\mu_1)}{M_1+pj+r} 
\notag \\[1.5mm]
& \hspace{4cm} \times
\qbins{M_2+L_2+((p'-p)j-r+s)/N-\frac{1}{2}(M_1+pj+r-\mu_{N-1})}{M_2-pj-r},
\notag 
\end{align}
with $(\vmu,\vet)$-systems
$$
\vmu+\vet=\frac{1}{2}({\mathcal I}\vmu+(M_1+pj)\ve_1+(M_2-pj)\ve_{N-1})
$$
for the first term of the right-hand side and
$$
\vmu+\vet=\frac{1}{2}({\mathcal I}\vmu+(M_1+pj+r)\ve_1+(M_2-pj-r)\ve_{N-1})
$$
for the second term of the right-hand side.
In Section~\ref{sec multinomials} we will show that in the limit 
when $M_1,M_2$ tend to infinity for fixed $M_{12}$ the above 
polynomials become proportional to the
one-dimensional configuration sums of solvable lattice models 
of Date et al.~\cite{DJKMO87,DJKMO88}, which are bounded
analogues of level-$N$ A$_1^{(1)}$ branching functions.

Using \eqref{sum}, it follows that
\begin{multline}\label{BurgetrafoN}
X_{r,r+Ns,\sigma}^{(p,p+Np'),N}(M_1,L_1,M_2,L_2) \\
=\sum_{i\in\Z} q^{i(i+M_{12})/N}
\qbin{L_1+L_2+M_2-i}{M_2-i} 
\sum_{\substack{\vn\in \Z^{N-1} \\
\frac{2i+M_{12}+\sigma N}{2N}+(C^{-1}\vn)_1\in \Z}}
q^{\vn C^{-1}\vn}\qbin{\vm+\vn}{\vn} \\
\times
X_{r,s}^{(p,p')}(i+M_{12},L_1-i+\frac{1}{2}m_1,i,L_2-M_{12}-i+\frac{1}{2}m_1),
\end{multline}
where on the right-hand side we assume the $(\vm,\vn)$-system
\begin{equation}\label{mni}
\vm+\vn=\frac{1}{2}(\I \vm+(2i+M_{12})\ve_1).
\end{equation}
Because of the conditions $L_1,L_2\geq 0$ in \eqref{sum}, a sufficiency
condition for the above transformation to hold is
\begin{equation}\label{suf}
\begin{split}
\left\lfloor\frac{L_1+M_{12}(N-1)/(2N)-s-r/N}{p'+p/N}\right\rfloor &\leq 
\left\lfloor\frac{L_1+M_{12}+r-s}{p'-p}\right\rfloor \\[1mm]
\left\lfloor\frac{L_2-M_{12}(N-1)/(2N)+s+r/N}{p'+p/N}\right\rfloor &\leq 
\left\lfloor\frac{L_2-M_{12}-r+s}{p'-p}\right\rfloor
\end{split}
\end{equation}
together with the inequalities obtained by setting $r=s=0$,
where we assumed that $p'>p$.
(The kernel of $X_{r,r+Ns,\sigma}^{(p,p+Np'),N}$ and 
of $X_{r,s}^{(p,p')}$ on either side of \eqref{BurgetrafoN}
is zero unless the summation variable $j$ lies 
in certain ranges. The above conditions make sure that in these
ranges of $j$ the conditions $L_1,L_2\geq 0$ of Theorem~\ref{thm ex burge}
apply).

Using the symmetry \eqref{symmetry} one also finds
\begin{multline}\label{trafo}
X_{s-M_{12},N(r+L_{12}+M_{12})+s-M_{12},\sigma}^{(p',Np+p'),N}
(M_1,L_1,M_2,L_2) \\
=\sum_{i\in\Z} q^{i(i+M_{12})/N}
\qbin{L_1+L_2+M_2-i}{M_2-i} 
\sum_{\substack{\vn\in \Z^{N-1} \\
\frac{2i+M_{12}+\sigma N}{2N}+(C^{-1}\vn)_1\in \Z}}
q^{\vn C^{-1}\vn}\qbin{\vm+\vn}{\vn} \\
\times X_{r,s}^{(p,p')}
(L_1-i+\frac{1}{2}m_1,i+M_{12},L_2-M_{12}-i+\frac{1}{2}m_1,i),
\end{multline}
where again \eqref{mni} holds. This time a sufficient condition is that
\begin{equation}\label{suf2}
\begin{split}
\left\lfloor\frac{L_2-M_{12}(N-1)/(2N)-s-r/N}{p'+p/N}\right\rfloor &\leq 
\left\lfloor\frac{L_1+M_{12}+r-s}{p'-p}\right\rfloor \\[1mm]
\left\lfloor\frac{L_1+M_{12}(N-1)/(2N)+s+r/N}{p'+p/N}\right\rfloor &\leq 
\left\lfloor\frac{L_2-M_{12}-r+s}{p'-p}\right\rfloor
\end{split}
\end{equation}
holds, as well as the inequalities obtained by setting $r=s=0$,
where again $p'>p$.

Again we consider the simpler case when $M_{12}=L_{12}=0$. Setting
$$X_{r,s,\sigma}^{(p,p'),N}(M,L,M,L)=X_{r,s,\sigma}^{(p,p'),N}(M,L),$$
the generalized Burge transformations \eqref{BurgetrafoN} and 
\eqref{trafo} simplify to
\begin{multline}\label{traf1}
X_{r,r+Ns,\sigma}^{(p,p+Np'),N}(M,L) \\
=\sum_{i=0}^M q^{i^2/N}\qbin{2L+M-i}{2L}
\sum_{\substack{\vn\in \Z^{N-1} \\
\frac{2i+\sigma N}{2N}+(C^{-1}\vn)_1\in \Z}}
\! q^{\vn C^{-1}\vn}\qbin{\vm+\vn}{\vn}
X_{r,s}^{(p,p')}(i,L-i+\frac{1}{2}m_1)
\end{multline}
and
\begin{multline}\label{traf2}
X_{s,Nr+s,\sigma}^{(p',Np+p'),N}(M,L) \\
=\sum_{i=0}^M q^{i^2/N}\qbin{2L+M-i}{2L} 
\sum_{\substack{\vn\in \Z^{N-1} \\
\frac{2i+\sigma N}{2N}+(C^{-1}\vn)_1\in \Z}}
\! q^{\vn C^{-1}\vn}\qbin{\vm+\vn}{\vn}
X_{r,s}^{(p,p')}(L-i+\frac{1}{2}m_1,i)
\end{multline}
both with ($\vm,\vn)$-system 
\begin{equation}\label{mni2}
\vm+\vn=\frac{1}{2}(\I \vm+2i\ve_1).
\end{equation}
The sufficiency conditions \eqref{suf} and \eqref{suf2}
(and their $r=s=0$ counterparts) reduce to the single condition
\begin{equation}\label{sufsym}
\left\lfloor\frac{L+s+r/N}{p'+p/N}\right\rfloor \leq 
\left\lfloor\frac{L-r+s}{p'-p}\right\rfloor.
\end{equation}

To end this section let us give some simple examples of our extensions to
the Burge transform, by finding the generalizations of equations
\eqref{nn}--\eqref{RR} to arbitrary $N$.
First, applying \eqref{traf1} to \eqref{initial} yields
\begin{equation*}
X_{0,N,\sigma}^{(1,2N+1),N}(M,L)=q^{L^2}
\sum_{\vm\in\Z^{N-1}} q^{\frac{1}{4}\vm T \vm} 
\qbin{L+M-\frac{1}{2}m_1}{2L}\qbin{\vm+\vn}{\vm},
\end{equation*}
with $\vm+\vn=\frac{1}{2}({\mathcal I}_T \vm+2L \ve_1)$ and
$({\mathcal I}_T)_{i,j}=\delta_{|i-j|,1}+\delta_{i,j}\delta_{i,1}$
the incidence matrix of the tadpole graph with $N-1$ nodes, and 
$T=2I-{\mathcal I}_T$ the corresponding Cartan-like matrix.
When $N$ is odd $\sigma=0$, $L\in\Z$ and $\vm\in 2\Z^{N-1}$.
When $N$ is even $m_{2i+1}\equiv 2L\equiv\sigma\pmod{2}$ and
$m_{2i}\equiv 0\pmod{2}$. The sufficiency condition \eqref{sufsym}
is satisfied.
Next applying \eqref{traf2} to \eqref{initial} yields 
\begin{equation*}
X_{1,1,\sigma}^{(2,N+2),N}(M,L)=\qbin{2L+M}{2L}\delta_{\sigma,0}
\end{equation*}
which, for $\sigma=0$, is a doubly bounded version of the Euler identity for
the level-$N$ string functions of type A$_1^{(1)}$.
Our third example follows after inserting \eqref{nn} into \eqref{traf2},
\begin{equation*}
X_{1,1,\sigma}^{(3,N+3),N}(M,L)=
\sum_{\vm\in \Z^N} q^{\frac{1}{4}\vm C \vm} 
\qbin{2L+M-\frac{1}{2}m_1}{2L}\qbin{\vm+\vn}{\vm},
\end{equation*}
with $(\vm,\vn)$-system $\vm+\vn=\frac{1}{2}(\I \vm+2L\ve_1)\in\Z^N$,
where $\I$ now is the incidence matrix of the A$_N$ Dynkin diagram.
When $N$ is odd $\sigma=0$, $L\in\Z$ and $\vm\in 2\Z^{N-1}$ and
when $N$ is even $m_{2i}\equiv 2L\equiv\sigma\pmod{2}$ and
$m_{2i+1}\equiv 0\pmod{2}$.
These identities are bounded analogues of identities for level-$N$
A$_1^{(1)}$ branching functions isomorphic to unitary minimal Virasoro 
characters.
Finally we use \eqref{traf1} and \eqref{Euler} to find
\begin{multline*}
X_{1,N+1,\sigma}^{(2,3N+2),N}(M,L)= \\
\sum_{i=0}^M q^{i^2/N} \qbin{2L+M-i}{2L}
\sum_{\substack{\vn\in \Z^{N-1} \\
\frac{2i+\sigma N}{2N}+(C^{-1}\vn)_1\in \Z}}
q^{\vn C^{-1} \vn}\qbin{\vm+\vn}{\vn}
\qbin{2L-i+m_1}{i}
\end{multline*}
where \eqref{mni2} holds.
As remarked before, for $N=1$ ($\sigma=0$) this is a doubly bounded 
version of the (first) Rogers--Ramanujan identity. For $N=2$ it becomes
\begin{equation*}
X_{1,3,\sigma}^{(2,8),2}(M,L)=
\sum_{i=0}^M \sum_{\substack{n=0 \\ n+i+\sigma~\text{even}}}^i  
q^{(i^2+n^2)/2} \qbin{2L+M-i}{2L}\qbin{i}{n}
\qbin{2L-n}{i}
\end{equation*}
which can be recognized as a doubly bounded version of
\begin{equation*}
\sum_{n=0}^{\infty} \frac{q^{n^2}(-q;q^2)_n}{(q^2;q^2)_n}=
\prod_{j=1}^{\infty}\frac{1}{(1-q^{8j-1})(1-q^{8j-4})(1-q^{8j-7})}
\end{equation*}
due to Slater~\cite{Slater52}
and related to the (first) G\"ollnitz--Gordon partition
identity~\cite{Goellnitz67,Gordon65},

\section{Special limits of Theorem~\ref{thm ex burge}}

\subsection{$q$-Multinomial coefficients}\label{sec multinomials}
In Refs.~\cite{Andrews94,Butler94,Kirillov95,Schilling96,Warnaar97} 
$q$-multinomial coefficients were introduced as $q$-analogues of the 
coefficients in the expansion 
\begin{equation*}
(1+x+x^2+\cdots+x^N)^L = 
 \sum_{a=-\frac{NL}{2}}^{\frac{NL}{2}} \bin{L}{a}_N x^{a+\frac{NL}{2}},
\end{equation*}
for $L\in\Zp$.
The $q$-multinomial coefficients are the generating function of a wide
class of combinatorial objects:
(i) unrestricted lattice paths related to the RSOS lattice models of Date
et al. with $H$-function statistic~\cite{DJKMO87,DJKMO88}, 
(ii) Durfee dissection partitions~\cite{Warnaar97}
and (iii) tabloids of shape $(N^L)$ and
content $(1^a 2^{NL-a})$ with the statistic ``value''~\cite{Butler94},
et cetera. 

Here we need the following explicit representation for the 
$q$-multinomials~\cite{Schilling96}
\begin{equation}\label{def T}
T_n^{(N)}(L,a)=\sum_{\substack{\vet\in\Z^{N-1}\\ 
\frac{L}{2}+\frac{a}{N}+(C^{-1}\vet)_1\in\Z}}
\frac{q^{\vet C^{-1}(\vet-\ve_n)}
(q)_L}{(q)_{\frac{L}{2}-\frac{a}{N}-(C^{-1}\vet)_1}
(q)_{\frac{L}{2}+\frac{a}{N}-(C^{-1}\vet)_{N-1}}(q)_{\vet}},
\end{equation}
where $L\in\Zp$, $2a\in\{-NL,-NL+2,\ldots,NL\}$ and $n\in\{0,1,\ldots,N-1\}$.
Repeated use of Newton's binomial expansion shows that $$
\lim_{q\to 1} T_n^{(N)}(L,a)=\bin{L}{a}_N$$
so that $T_n^{(N)}(L,a)$ is indeed a $q$-analogue of the multinomial 
coefficient.

Theorem~\ref{thm ex burge} provides a new representation of the 
$q$-multinomials when $n=0$. To see this we
let $M$ tend to infinity in \eqref{sum} resulting in
\begin{multline*}
\sum_{i=0}^{\infty} q^{i(i+\ell)/N}  
\sum_{\substack{\vn\in \Z^{N-1} \\ 
\frac{2i+\ell+\sigma N}{2N}+(C^{-1} \vn)_1\in\Z}}
q^{\vn C^{-1} \vn} \qbins{\vm+\vn}{\vn} 
\qbins{L_1+\frac{1}{2}m_1}{i+\ell} \qbins{L_2+\frac{1}{2}m_1}{i} \\
=\sum_{\substack{\vet\in \Z^{N-1} \\ 
\frac{\ell+\sigma N}{2N}+(C^{-1}\vet)_1\in\Z}}
\frac{q^{\vet C^{-1}\vet} (q)_{L_1+L_2}}{
(q)_{L_1-\frac{\ell}{2N}-\frac{\ell}{2}-(C^{-1}\vet)_1} 
(q)_{L_2+\frac{\ell}{2N}+\frac{\ell}{2}-(C^{-1}\vet)_{N-1}}
(q)_{\vet}}.
\end{multline*}

If we now set $L_1=\frac{1}{2}(L+\ell)$ and $L_2=\frac{1}{2}(L-\ell)$
(so that $\sigma\equiv L\pmod{2}$) and compare with the
right-hand side of \eqref{def T}, we find that
\begin{multline}\label{Tnew}
T_0^{(N)}(L,\ell/2)=\\
\sum_{i=0}^{\infty}  q^{i(i+\ell)/N}  
\sum_{\substack{\vn\in\Z^{N-1} \\ 
\frac{L}{2}+\frac{2i+\ell}{2N}+(C^{-1} \vn)_1\in\Z}}
q^{\vn C^{-1} \vn} \qbins{\vm+\vn}{\vn} 
\qbins{\frac{1}{2}(L+\ell+m_1)}{i+\ell}
\qbins{\frac{1}{2}(L-\ell+m_1)}{i},
\end{multline}
with $\vm$ given by \eqref{mn}.

When $N=1$ the above decomposition of the $q$-multinomial coefficients
reduces to the $q$-Chu--Vandermonde sum \eqref{qCV} and a combinatorial
interpretation is easily given as follows.
The $q$-binomial $\qbins{m+n}{m}$ is the generating function
of partitions that fit in a box of dimension $m$ times $n$.
Hence the summand on the left-hand side of \eqref{qCV}
is the generating function of
partitions that fit in a box of dimension $L_1-\ell$ times $L_2+\ell$ which
have a Durfee rectangle of size $i$ by $i+\ell$ 
(maximal rectangle of the Ferrers graph that has a
horizontal excess of $\ell$ nodes). Summing over $i$ removes the
Durfee rectangle restriction resulting in the right-hand side.
It seems an interesting problem to also explain the $q$-multinomial
decomposition \eqref{Tnew} combinatorially.

There is a corresponding formula for $1\leq n<N-1$ which, however,
is less appealing (and which we will not prove here)
\begin{multline*}
T_n^{(N)}(L,(n-\ell)/2)-q^{(\ell+1)/N}T_n^{(N)}(L,(n+\ell+2)/2) \\
= \sum_{i=0}^{\infty} q^{i(i+\ell)/N}
\hspace{-2mm}
\sum_{\substack{n\in\Z^{N-1}\\
\frac{L}{2}+\frac{2i+\ell-n}{2N}+(C^{-1} \vn)_1 \in \Z}}
\hspace{-5mm}
q^{\vn C^{-1} (\vn-\ve_{N-n})} \qbins{\vm+\vn}{\vn} \\[2mm]
\times \Bigl(
\qbins{\frac{1}{2}(L+\ell+m_1}{i+\ell}
\qbins{\frac{1}{2}(L-\ell+m_1)}{i}
-\qbins{\frac{1}{2}(L+\ell+2+m_1)}{i+\ell+1}
\qbins{\frac{1}{2}(L-\ell-2+m_1)}{i-1} \Bigr)
\end{multline*}
with
\begin{equation*}
\vm+\vn=\frac{1}{2}(\I \vm+(2i+\ell)\ve_1+\ve_{N-n}).
\end{equation*}
Although this identity has the structure
$f(L,\ell)-q^{(\ell+1)/N}f(L,-\ell-2)=
g(L,\ell)-q^{(\ell+1)/N}g(L,-\ell-2)$, it is not true that
$f(L,\ell)=g(L,\ell)$.

To conclude our discussion of the $q$-multinomial coefficients, let
us point out that the polynomials defined in equation \eqref{XN} are
related to one-dimensional configuration sums of lattice models of
Date et al~\cite{DJKMO87,DJKMO88}. Let $L\in\Z$ and choose 
$$L_1=\frac{1}{2}\Bigl(L-M_{12}-\frac{r-s}{N}\Bigr) \qquad
L_2=\frac{1}{2}\Bigl(L+M_{12}+\frac{r-s}{N}\Bigr)$$
so that $\sigma=0,1$ is fixed by the condition that
$L-(r-s)/N+\sigma$ is even. Then
\begin{align*}
\lim_{\substack{M_1,M_2\to\infty \\ \text{$M_{12}$ fixed}}}(q)_{2L}&
X_{r,s,\sigma}^{(p,p'),N}(M_1,L_1,M_2,L_2) \\
=\sum_{j=-\infty}^{\infty}\Bigl\{&
q^{\frac{j}{N}(pp'j+p'(M_{12}+r)-ps)} 
T_0^{(N)}\bigl(L,\frac{1}{2}(r+M_{12}-s)+p'j\bigr)\\
-&q^{\frac{1}{N}(pj+M_{12}+r)(p'j+s)}
T_0^{(N)}\bigl(L,\frac{1}{2}(r+M_{12}+s)+p'j\bigr)\Bigr\},
\end{align*}
which, for $p'=p+N$, is proportional to the configuration sums of the
models of Date et al. in the representation obtained in 
\cite[Eq. (3.15)]{Schilling96}.

\subsection{Bailey's lemma}\label{sec bailey}
In this section we show that the limit $L_1,L_2\to\infty$
of Theorem \ref{thm ex burge} gives rise to the higher-level 
Bailey lemma (or more precisely the higher-level conjugate 
Bailey pairs) of Refs.~\cite{SW97,SW98}.

Bailey's lemma~\cite{Bailey49} is an elegant tool to prove
$q$-series identities such as the famous Rogers--Ramanujan identities.
Let $\alpha=\{\alpha_L\}_{L\geq 0},\beta=\{\beta_L\}_{L\geq 0}$
be a pair of sequences that satisfies
\begin{equation*}
\beta_L=\sum_{i=0}^L \frac{\alpha_i}{(q)_{L-i}(aq)_{L+i}}.
\end{equation*}
Such a pair is called a 
Bailey pair relative to $a$. Recalling the definition \eqref{cBp}
of a conjugate Bailey pair, it follows by a simple interchange of sums that
\begin{equation}\label{abcd}
\sum_{L=0}^{\infty} \alpha_L \gamma_L
=\sum_{L=0}^{\infty} \beta_L \delta_L.
\end{equation}
Many known $q$-series identities follow from \eqref{abcd}
after substitution of suitable Bailey and conjugate Bailey pairs.

Now let $L_1,L_2$ tend to infinity in \eqref{sum} and replace
$i\to i-L$, $\ell\to\ell+2L$ and $M\to M-L$. This yields
\begin{multline*}
\sum_{i=L}^M \frac{q^{i(i+\ell)/N}}{(q)_{i-L}(q)_{i+L+\ell}(q)_{M-i}}
\sum_{\substack{\vn\in\Z^{N-1}\\ 
 \frac{2i+\ell+\sigma N}{2N}+(C^{-1}\vn)_1\in\Z}}
q^{\vn C^{-1}\vn}\qbin{\vm+\vn}{\vn}\\
=\frac{q^{L(L+\ell)/N}}{(q)_{M-L}(q)_{M+L+\ell}}
\sum_{\substack{\vet\in\Z^{N-1}\\ 
 \frac{2L+\ell+\sigma N}{2N}+(C^{-1}\vet)_1\in\Z}}
q^{\vet C^{-1}\vet}\qbin{\vmu+\vet}{\vet},
\end{multline*}
with $(\vm,\vn)$-system \eqref{mn} and $(\vmu,\vet)$-system
\begin{equation}\label{muet1}
\vmu+\vet=\frac{1}{2}(\I\vmu+(M+L+\ell)\ve_1+(M-L)\ve_{N-1}).
\end{equation}
Comparing with \eqref{cBp} one reads off the following conjugate Bailey pair
(which is the special case $\lambda=0$ of \cite[Corollary 2.1]{SW98})
\begin{equation*}
\begin{split}
\gamma_L&=\frac{a^{L/N}q^{L^2/N}}{(q)_{M-L}(aq)_{M+L}}
\sum_{\substack{\vet\in\Z^{N-1}\\ \frac{2L+\ell+\sigma N}{2N}
+(C^{-1}\vet)_1\in\Z}}
q^{\vet C^{-1}\vet}\qbin{\vmu+\vet}{\vet},\\
\delta_L&=\frac{a^{L/N}q^{L^2/N}}{(q)_{M-L}}
\sum_{\substack{\vn\in\Z^{N-1}\\ 
\frac{2L+\ell+\sigma N}{2N}+(C^{-1}\vn)_1\in\Z}}
q^{\vn C^{-1}\vn}\qbin{\vm+\vn}{\vn},
\end{split}
\end{equation*}
with $a=q^{\ell}$ and where \eqref{muet1} and
$\vm+\vn=\frac{1}{2}(\I\vm+(2L+\ell)\ve_1)$ hold.

\subsection{String functions}\label{sec string}
Taking the limit $L_1,L_2,M\to\infty$ in Theorem~\ref{thm ex burge}
we obtain
\begin{multline}\label{lim LM}
\sum_{i=0}^{\infty}\frac{q^{i(i+\ell)/N}}{(q)_i(q)_{i+\ell}}
\sum_{\substack{\vn\in\Z^{N-1}\\
\frac{2i+\ell+\sigma N}{2N}+(C^{-1}\vn)_1\in\Z}}
q^{\vn C^{-1}\vn}\qbin{\vm+\vn}{\vn}\\
=\frac{1}{(q)_{\infty}}
\sum_{\substack{\vet\in\Z^{N-1}\\ \frac{\ell+\sigma N}{2N}
+(C^{-1}\vet)_1\in\Z}}
\frac{q^{\vet C^{-1}\vet}}{(q)_{\vet}}.
\end{multline}
It was shown in Refs.~\cite{ANOT96,BLS95,NY96a,NY96b,SW99} that
the left-hand side is proportional to a 
level-$N$, A$_1^{(1)}$ string function $C_{m,\ell}^N$ defined as follows.
Let
\begin{equation*}
\Theta_{n,m}(z,q)=\sum_{j\in\Z+n/2m} q^{mj^2}z^{-mj}
\end{equation*}
be the classical theta function of degree $m$ and characteristic $n$.
The A$_1^{(1)}$ character of the highest weight module of highest
weight $(N-\ell)\Lambda_0+\ell\Lambda_1$ (where $\Lambda_0$ and 
$\Lambda_1$ are the fundamental weights of A$_1^{(1)}$ and
$0\le \ell\le N$) is given by
\begin{equation*}
\chi_{\ell}(z,q)=\frac{\sum_{\sigma=\pm 1}\sigma
\Theta_{\sigma(\ell+1),N+2}(z,q)}
{\sum_{\sigma=\pm 1}\sigma\Theta_{\sigma,2}(z,q)}.
\end{equation*}
The level-$N$ A$_1^{(1)}$ string functions are defined by the expansion
\begin{equation*}
\chi_{\ell}(z,q)=\sum_{m\in 2\Z+\ell} C_{m,\ell}^N(q) 
q^{\frac{m^2}{4N}} z^{-\frac{1}{2}m}.
\end{equation*}
According to the above-cited references
\begin{align*}
C^N_{m,\ell}(q)
&=q^{\frac{(\ell+1)^2}{4(N+2)}-\frac{m^2}{4N}-\frac{1}{8}}
\sum_{i=0}^{\infty}\frac{X_{\ell+1}^{N+2}(2i+m)}{(q)_i(q)_{i+m}} \\
&=q^{\frac{(\ell+1)^2}{4(N+2)}-\frac{\ell^2}{4N}-\frac{1}{8}}
\sum_{i=0}^{\infty}\frac{q^{i(i+m)/N}}{(q)_i(q)_{i+m}}
\sum_{\substack{\vn\in\Z^{N-1}\\
\frac{2i+m+\ell}{2N}+(C^{-1}\vn)_1\in\Z}}
q^{\vn C^{-1}(\vn-\ve_{\ell})}\qbin{\vm+\vn}{\vn},
\end{align*}
with $\vm+\vn=\frac{1}{2}(\I \vm+(2i+m)\ve_1+\ve_{\ell})$ and 
$X_s^p(L)$ a one-dimensional configuration sum
of the $(p-1)$-state Andrews--Baxter--Forrester model in regime I,
\begin{equation*}
X_s^p(L)= \sum_{j=-\infty}^{\infty} q^{j(pj+s)}
\biggl\{\qbin{L}{\frac{1}{2}(L-s+1)-pj}
-\qbin{L}{\frac{1}{2}(L-s-1)-pj}\biggr\}.
\end{equation*}

Comparing with \eqref{lim LM} we obtain the following
expression of the string function
\begin{equation*}
C_{m,\sigma N}^N(q)=\frac{q^{\frac{1}{4(N+2)}-\frac{1}{8}}}{(q)_\infty}
\sum_{\substack{\vet\in\Z^{N-1}\\ \frac{m+\sigma N}{2N}+(C^{-1}\vet)_1\in\Z}} 
\frac{q^{\vet C^{-1}\vet}}{(q)_{\vet}},
\end{equation*}
which was first derived by Lepowsky and Primc~\cite{LP85}.

\subsection*{Acknowledgements}
We thank Omar Foda and Trevor Welsh for discussions on the Burge transform.
The first author was supported by the ``Stichting Fundamenteel
Onderzoek der Materie''. The second author was
supported by a fellowship of the Royal
Netherlands Academy of Arts and Sciences.

\bibliographystyle{amsplain}

\begin{thebibliography}{99}

\bibitem{Andrews76}
G. E. Andrews,
\textit{The Theory of Partitions},
Encyclopedia of Mathematics and its Applications, Vol.~2,
(Addison-Wesley, Reading, Massachusetts, 1976).

\bibitem{Andrews94}
G. E. Andrews,
\textit{Schur's theorem, Capparelli's conjecture and $q$-trinomial
coefficients},
Contemp. Math. \textbf{166} (1994) 141--154.

\bibitem{AB87}
G. E. Andrews and R. J. Baxter,
\textit{Lattice gas generalization of the hard
hexagon model. III. $q$-Trinomial coefficients},
J. Stat. Phys. \textbf{47} (1987) 297--330.

\bibitem{ANOT96}
T. Arakawa, T. Nakanishi, K. Oshima and A. Tsuchiya,
\textit{Spectral decomposition of path space in solvable lattice model},
Comm. Math. Phys. \textbf{181} (1996) 157--182.

\bibitem{Bailey49}
W. N. Bailey,
\textit{Identities of the Rogers--Ramanujan type},
Proc. London Math. Soc. (2) \textbf{50} (1949) 1--10.

\bibitem{BLS95}
P. Bouwknegt, A. W. W. Ludwig and K. Schoutens,
\textit{Spinon basis for higher level $SU(2)$ WZW models},
Phys. Lett. B \textbf{359} (1995) 304--312.

\bibitem{Burge93}
W.~H.~Burge,
\textit{Restricted partition pairs},
J. Combin. Theory Ser. A \textbf{63} (1993) 210--222.

\bibitem{Butler94}
L.~M.~Butler,
{\em Subgroup lattices and symmetric functions},
Memoirs of the Amer. Math. Soc., no. 539, vol. {\bf 112} (1994).

\bibitem{Carlitz74}
L.~Carlitz, 
\textit{Remark on a combinatorial identity},
J. Combin. Theory Ser. A \textbf{17} (1974) 256--257. 

\bibitem{DJKMO87}
E. Date, M. Jimbo, A. Kuniba, T. Miwa and M. Okado,
\textit{Exactly solvable SOS models: local height probabilities and theta 
function identities},
Nucl. Phys. B \textbf{290} (1987) 231--273.

\bibitem{DJKMO88}
E. Date, M. Jimbo, A. Kuniba, T. Miwa and M. Okado,
\textit{Exactly solvable SOS models. II. Proof of the star-triangle 
relation and combinatorial identities},
Adv. Stud. Pure Math. \textbf{16} 17--122.

\bibitem{FLW98}
O. Foda, K. S. M. Lee and T. A. Welsh,
\textit{A Burge tree of Virasoro-type polynomial identities},
Int. J. Mod. Phys. A \textbf{13} (1998) 4967--5012.

\bibitem{GR90}
G.~Gasper and M.~Rahman,
\textit{Basic Hypergeometric Series},
Encyclopedia of Mathematics and its Applications, Vol.~35,
(Cambridge University Press, Cambridge, 1990).

\bibitem{GK97}
I.~M.~Gessel and C.~Krattenthaler,
\textit{Cylindric partitions},
Trans. Amer. Math. Soc. \textbf{349} (1997) 429--479.

\bibitem{Goellnitz67}
H.~G\"ollnitz,
\textit{Partitionen mit Differenzenbedingungen},
J. Reine Angew. Math. \textbf{225} (1967) 154--190.
 
\bibitem{Gordon65}
B. Gordon,
\textit{Some continued fractions of the Rogers--Ramanujan type},
Duke Math. J. \textbf{31} (1965) 741--748.

\bibitem{Gould72}
H.~W.~Gould,
\textit{A new symmetrical combinatorial identity},
J. Combin. Theory Ser. A \textbf{13} (1972) 278--286. 

\bibitem{Kirillov95}
A. N. Kirillov,
\textit{Dilogarithm identities},
Prog. Theor. Phys. Suppl. \textbf{118} (1995) 61--142.

\bibitem{LP85}
J. Lepowsky and M. Primc,
\textit{Structure of the standard modules for the affine Lie algebra
A$_1^{(1)}$}, 
Contemp. Math. Vol. 46 (AMS, Providence, 1985).

\bibitem{NY96a}
A. Nakayashiki and Y. Yamada,
\textit{Crystallizing the spinon basis},
Comm. Math. Phys. \textbf{178} (1996) 179--200.
 
\bibitem{NY96b}
A. Nakayashiki and Y. Yamada,
\textit{Crystalline spinon basis for RSOS models},
Int. J. Mod. Phys. A \textbf{11} (1996) 395--408.

\bibitem{Schilling96}
A. Schilling,
\textit{Multinomials and polynomial bosonic forms for the 
branching functions of the $\widehat{su}(2)_M\times
\widehat{su}(2)_N/\widehat{su}(2)_{M+N}$ conformal coset models},
Nucl. Phys. B \textbf{467} (1996) 247--271.

\bibitem{SW97}
A. Schilling and S. O. Warnaar,
\textit{A higher-level Bailey lemma},
Int. J. Mod. Phys. B \textbf{11} (1997) 189--195.

\bibitem{SW98}
A. Schilling and S. O. Warnaar,
\textit{A higher level Bailey lemma: Proof and application},
The Ramanujan Journal \textbf{2} (1998) 327--349.

\bibitem{SW99}
A. Schilling and S. O. Warnaar,
\textit{Conjugate Bailey pairs. From configuration sums and
fractional-level string functions to Bailey's lemma},
preprint math.QA/9906092.

\bibitem{Slater52}
L. J. Slater,
\textit{Further identities of the Rogers--Ramanujan type},
Proc. London Math. Soc. (2) \textbf{54} (1952) 147--167.

\bibitem{Warnaar97}
S. O. Warnaar,
\textit{The Andrews-Gordon identities and $q$-multinomial coefficients},
Comm. Math. Phys. \textbf{184} (1997) 203--232.

\end{thebibliography}

\end{document}